\documentclass[12pt]{article}
\usepackage{subeqn}
\usepackage{graphicx}
\usepackage{psfrag}

\advance \topmargin by -\headheight
\advance \topmargin by -\headsep
     
\evensidemargin \oddsidemargin
\def\pr{\prime}
\newcommand\jbk{Joseph B. Keller
\\{\small  Departments of Mathematics and Mechanical Engineering}\\
{\small Stanford University}\\{\small  Stanford, CA  94305-2125}
}

\newcommand\be{\begin{equation}}
\newcommand\ee{\end{equation}}
\newcommand\bea{\begin{eqnarray}}
\newcommand\eea{\end{eqnarray}}
\newcommand\nn{\nonumber}

\newcommand\sce{\setcounter{equation}{0}}
\newcommand\lam{\lambda}
\def\pmb#1{\setbox0=\hbox{#1}%
  \kern-.025em\copy0\kern-\wd0
   \kern.05em\copy0\kern-\wd0
   \kern-.025em\raise.0433em\box0 }

\begin{document}
\title{\bf  Multiple eigenvalues   }

\author{\jbk}
\date{June 24, 2007}
\maketitle

\begin{abstract}
The dimensions of sets of matrices of various types, with specified eigenvalue multiplicities, are determined.  The dimensions of the sets of matrices with given Jordan form and with given singular value multiplicities are also  found.  Each corresponding codimension is the number of conditions which a matrix of the given type must satisfy in order to have the specified multiplicities.
\end{abstract}

\section{Introduction}

We shall determine the dimensions of the sets of diagonalizable, normal, Hermitian, unitary, and real symmetric matrices having one eigenvalue of multiplicity $k_1$, another of multiplicity $k_2$, etc.  We shall also find the dimensions of the set of  square matrices having a specified Jordan form, and of the set of real rectangular matrices having largest singular  value of multiplicity $k_1$,  next largest of multiplicity $k_2$, etc.  In each case the codimension of the specified set gives the number of conditions which the elements of the matrix must satisfy in order to be in that set.  For example, the elements of a real symmetric matrix must satisfy two conditions  in order to have a  double eigenvalue, which is a result of von Neuman and Wigner \cite{vN&W}.

To determine these dimensions, we first find the set of all matrices which transform a given matrix into its  specified diagonal form, Jordan canonical form, or singular value form.  Then we calculate the dimensions of these sets of transforming matrices, and use them to obtain the desired results.

Our main results are listed in Table 2.  The dimensions of some common sets of matrices are given in Table 1.

This work was stimulated by a lecture of Beresford Parlett, based upon \cite{BP}, which is related to the work of Lax \cite{Lax}.  Overton and Womersley \cite{O&W} obtained results like some of those in the present paper.


\section{Diagonalizable Matrices in $C^{nn}$}
\sce
Let $C^{nn}$ be the set of complex   square matrices of order $n$.  In $C^{nn}$ a  matrix $A$  is similar to $J$ if and only if there is an invertible $T$ such that $A=TJT^{-1}$.  If in addition  $A= T_1 JT^{-1}_1$ then $TJT^{-1} = T_1 JT^{-1}_1$, so $T^{-1}_1 TJ=JT^{-1}_1 T$.  Thus $S=T^{-1}_1 T$ commutes with $J$:
\be
SJ=JS.
\label{2.1}
\ee
Since $T_1 =TS^{-1}$, we have the following result:

\noindent {\bf Lemma 1.}  If $A=TJT^{-1}$ then $A=T_1 JT^{-1}_1$ if and only if $T_1 =TS^{-1}$ where $S$ is invertible and satisfies (\ref{2.1}).

By using Lemma 1, we can calculate the dimension of the set of $A$ which are similar to $J$.  This dimension is just the dimension of the set of invertible $T$ minus the dimension of the set of invertible $S$ satisfying (\ref{2.1}).  Thus we have

\noindent {\bf Lemma 2.}
\be
\mbox{ dim }\left\{ A\left| \right. A=TJT^{-1} \right\} =\mbox{ dim } \left\{ T\left|\right.  T T^{-1}=I\right\} -\mbox{ dim } \left\{ S\left|\right. S S^{-1}=I, \quad  SJ = JS\right\}.
\label{2.2}
\ee

\noindent The complex dimension of the set of invertible $T$ in $C^{nn}$ is $n^2$, and its real dimension is $2n^2$.

We shall now find the set of  $S$ which commute with   $J$ when $J=\Lambda$ is diagonal.  We assume that $C^{nn}$ has $I$ distinct eigenvalues $\lambda_1, \lambda_2,\cdots,\lambda_I$ with multiplicities $k_1, k_2,\cdots,k_I$, where $k_1+ k_2 +\cdots+k_I =n$.  We write $\Lambda$ in block diagonal form
\be
\Lambda=\left( \lambda_1 I_1, \lambda_2 I_2,\cdots,\lambda_I I_I\right).
\label{2.3}
\ee
Here $I_i$ is the identity matrix of order $k_i$.  We partition $S$ into blocks $S_{ij}$,  of order $k_i$ by $k_j$, so that $S$ is conformable with $\Lambda$.   By setting $J=\Lambda$, and using (\ref{2.3}) and  this partition of $S$ in (\ref{2.1}), we get 
\be
S_{ij} \lambda_j =\lambda_i S_{ij} \qquad \mbox{ (no summation) }.
\label{2.4}
\ee
Since the $\lambda_i$ are distinct, it follows that $S_{ij}=0$ for $i\not=j$.  Thus $S$ is block diagonal with the block $S_{ii}$ square of order $k_i$, and it is unrestricted by (\ref{2.4}).  Therefore the complex dimension  of the set of $S_{ii}$ is $k^2_i$, and the complex dimension of the set of $S$ is the sum of the $k^2_i$.  We summarize this result as follows:

\noindent{\bf Lemma 3.} 
\be 
\mbox{ complex dim } \left\{ S\left|\right. SS^{-1} =I, S\Lambda =\Lambda S, \quad \mbox{mult } \Lambda =\left( k_1, k_2,\cdots,k_I\right)\right\} =\sum\limits^I_{i=1} k^2_i.
\label{2.5}
\ee

Now we use (\ref{2.5}) in (\ref{2.2}) to get 

\noindent{\bf Theorem 1.}
\be
\mbox{ complex dim } \left\{ A\left| \right. A=T\Lambda T^{-1}, \quad \mbox{mult } \Lambda =\left( k_1, \cdots,k_I\right)\right\} =n^2 -\sum\limits^I_{i=1} k^2_i.
\label{2.6}
\ee
This theorem gives the dimension of the set of diagonalizable $A$ with specified eigenvalues and eigenvalue multiplicities.  To obtain the dimension of the set of diagonalizable $A$ with specified multiplicities, we add to the right side of (\ref{2.6}) the dimension $I$ of the set of eigenvalues:

\noindent {\bf Corollary 1.1.}
\be
\mbox{ complex dim } \left\{ A\left|\right. A \mbox{ diagonalizable}, \quad \mbox{mult}\left( k_1,\cdots, k_I\right) \right\} =n^2 - \sum\limits^I_{i=1} k^2_i +I =n^2- \sum\limits^I_{i=1} \left( k^2_i -1\right).
\label{2.7}
\ee
The complex codimension of the set of diagonalizable $A$ with the specified eigenvalue multiplicities is just $n^2$, the dimension of $C^{nn}$, minus the dimension given in  in (\ref{2.7}):

\noindent { \bf Corollary 1.2.}
\be
\mbox{ complex codim } \left\{ A\left|\right. A \mbox{ diagonalizable}, \quad \mbox{mult}\left( k_1,\cdots,k_I\right)\right\} =\sum\limits^I_{i=1} \left( k^2_i -1\right).
\label{2.8}
\ee

As an application of Corollary 1.2, we consider diagonalizable $A$ with one $k$-fold eigenvalue and all others simple.  Then $I=n-k+1$, $ k_1 =k$, $k_2 =\cdots=k_{n-k+1}=1$ and (\ref{2.8}) yields
\be
\mbox{ complex codim } \left\{ A\left|\right. A \mbox{ diagonalizable}, \quad \mbox{mult }=\left(k, 1,\cdots, 1\right) \right\} =k^2 -1.
\label{2.9}
\ee

\section{Normal Matrices in $C^{nn}$}
\sce

A matrix $A$ in $C^{nn}$ is normal if it commutes with its adjoint $A^\ast$ : $AA^\ast =A^\ast A$.  Every normal matrix is similar to a diagonal matrix $\Lambda$, $A = U\Lambda U^\ast$, where $U$ is unitary $\left( U^{-1} =U^\ast\right)$.  Therefore $A$ is diagonalizable, so Lemmas 1 and 2 apply with the additional condition that $S$ is unitary.  This condition follows from the  definition $S=T^{-1}_1 T =U^{-1}_1 U$, since both $T=U$ and $T_1=U_1$ are unitary.  Equation (\ref{2.4}) and the consequence that $S$ is block diagonal still apply, but now the block $S_{ii}$ must be unitary.  The real dimension of the set of unitary matrices of order $k_i$ is $k^2_i$ so the real dimension of the set of unitary $S$ that commute with $\Lambda$ is the sum of the $k^2_i$.  Thus instead of Lemma 3 we have

\noindent {\bf Lemma 4.}
\be
\mbox{ real dim } \left\{ S\left| \right. SS^\ast =I, \quad S\Lambda =\Lambda S, \quad \mbox{mult }\Lambda =\left( k_1,k_2,\cdots,k_I\right)\right\} =\sum\limits^I_{i=1} k^2_i.
\label{3.1}
\ee

We now use in (\ref{2.2}) both (\ref{3.1}) and the fact that $T= U$ is unitary.  Since the real dimension of the set of $U$ is $n^2$, we obtain

\noindent {\bf Theorem 2.}
\be
\mbox{ real dim } \left\{ A\left| \right. A=U\Lambda U^\ast, \quad \mbox{mult }\Lambda =\left( k_1,\cdots, k_I\right)\right\} =n^2-\sum\limits^I_{i=1} k^2_i.
\label{3.2}
\ee
Upon adding to (\ref{3.2}) the real dimension of the set of $\Lambda$, which is $2I$, we get the  dimension of the set of normal $A$ in $C^{nn}$ with   specified multiplicities:

\noindent{\bf Corollary 2.1.}
\be
\mbox{ real dim } \left\{A\left| \right. A \mbox{ normal}, \quad \mbox{mult}\left( k_1,\cdots,k_I\right)\right\}=n^2 -\sum\limits^I_{i=1} k^2_i + 2I.
\label{3.3}
\ee

The dimension of the set of normal matrices is the maximum value of the right side of  (\ref{3.3}).  This is achieved when $I=n$ and each $k_i =1$, which gives
\be
\mbox{ real dim } \left\{ A\left|\right. A \mbox{ normal }\right\} =n^2 + n.
\label{3.4}
\ee
Upon subtracting the dimension in (\ref{3.3}) from the real dimension (\ref{3.4}) of the set of normal matrices, we get

\noindent{\bf Corollary 2.2.}
\be
\mbox{ real  codim }\left\{ A\left| \right. A \mbox{ normal}, \quad \mbox{mult}\left(
 k_1,\cdots, k_I\right)\right\} =
n+\sum\limits^I_{i=1} k^2_i -2I =\sum\limits^I_{i=1} \left( k_i-1\right) \left( k_i+2\right).
\label{3.5}
\ee

When $A$ is normal, with one $k$-fold eigenvalue and all the others simple, (\ref{3.5}) yields 
\be
\mbox{ real codim } \left\{ A\left|\right. A \mbox{ normal}, \quad \mbox{mult}\left( k, 1,\cdots, 1\right)\right\} = \left( k-1\right) \left( k+2\right).
\label{3.6}
\ee

\section{Hermitian matrices}
\sce

A matrix $A$ in $C^{nn}$ is Hermitian if it equals its adjoint, $A=A^\ast$, so $A$ is also normal.  
Therefore the results of section 3 up to and including Theorem 2 apply to Hermitian matrices.  The eigenvalues of an Hermitian matrix are real, so the dimension of the set of $\lambda_i$ is $I$.  Upon adding $I$ to (\ref{3.2}) we get the dimension of the set of Hermitian $A$ with multiplicities $k_1,\cdots,k_I$:

\noindent{\bf Corollary 2.3.}
\be
\mbox{real dim} \left\{ A\left| \right.   A \mbox{ Hermitian}, \,\mbox{mult}\left( k_1,k_2,\cdots,k_I\right) \right\} =n^2 - \sum\limits^I_{i=1} k^2_i +I =n^2 -\sum\limits^I_{i=1} \left( k^2_i -1\right).
\label{4.1}
\ee
This result was given by von Neuman and Wigner \cite{vN&W}.

Upon subtracting the dimension (\ref{4.1}) from $n^2$, the dimension of the set of Hermitian matrices of order $n$, we get

\noindent{\bf Theorem 3.}  In the set of Hermitian matrices of order $n$, the real codimension of the subset of those having $I$ distinct eigenvalues with multiplicities $k_1,\cdots,k_I$, where $k_1+\cdots+k_I=n$, is 
\be
\mbox{real codim} \left\{ A\left| \right. A=A^\ast , \,\mbox{mult}\left(k_1,\cdots, k_I\right) \right\} =\sum\limits^I_{i=1} \left( k^2_i-1\right).
\label{4.2}
\ee
 
When $k_1=k$ and the right side of (\ref{4.2}) is minimized over the other $k_i$, the result is $k^2-1$.  This is the number  of  real conditions on the elements of an Hermitian matrix for it to have a $k$-fold eigenvalue:

\noindent{\bf Corollary 3.1.}
\be
\mbox{ real codim }\left\{ A\left| \right. A\mbox{ Hermitian}, \,\mbox{mult}\left( k, 1,\cdots,1\right)\right\} =k^2-1.
\label{4.3}
\ee

A matrix $A$ in $C^{nn}$ is skew-Hermitian if $A=-A^\ast$.  Then $iA$ is Hermitian, so the results of this section yield corresponding results for skew-Hermitian matrices.

\section{Unitary matrices}
\sce

A matrix $A$ in $C^{nn}$ is   unitary if $AA^\ast =I$, so $A$ is also normal.  Therefore the results of section~3 through Theorem 2 apply to unitary    matrices.  Since the eigenvalues of a unitary matrix have absolute   value one, the dimension of the   set of $\lambda_i$ is $I$.  This is the same  as the dimension of the set of $\lambda_i$ for Hermitian matrices, so the dimension of the set of unitary $A$ with multiplicities $k_1,\cdots,k_I$ is also given by the right side of (\ref{4.1}):

\noindent{\bf Corollary 2.4.}
\be
\mbox{real dim} \left\{ A\left| \right. A\mbox{ unitary}, \,\mbox{mult}\left(k_1,\cdots,k_I\right)\right\} =n^2 -\sum\limits^I_{i=1} \left( k^2_i-1\right).
\label{5.1}
\ee

When all $k_i=1$, (\ref{5.1}) yields $n^2$, which is the dimension of the set of unitary matrices.  Therefore the results for unitary matrices, analogous to (\ref{4.2}) and (\ref{4.3}), are 

\noindent{\bf Theorem 4.}  In the set of unitary matrices of order $n$, the real codimension of the subset of those having $I$ distinct eigenvalues with multiplicities $k_1,\cdots,k_I$, where $k_1+\cdots+k_I=n$, is
\be
\mbox{real codim} \left\{ A\left| \right. A \mbox{ unitary}, \,\mbox{mult}  \left( k_1,\cdots,k_I\right)\right\}=\sum\limits^I_{i=1} \left( k^2_i -1\right).
\label{5.2}
\ee
The number of conditions on the elements of a unitary matrix for it to have a $k$-fold eigenvalue, is given by

\noindent{\bf Corollary 4.1.}
\be
\mbox{ real codim }\left\{ A\left| \right. A \mbox{ unitary}, \, \mbox{  mult }\left( k, 1,\cdots,1\right)\right\} 
=k^2-1.
\label{5.3}
\ee

\section{Real symmetric matrices}
\sce
 A matrix $A$ is real symmetric if $A$ is real and $A=A^T$, where $A^T$ is the transpose of $A$.  Therefore it is Hermitian, normal, and diagonalizable.  Every real symmetric $A$ is similar to a real diagonal matrix $\Lambda$, $A=O\Lambda O^T$, where $O$ is orthogonal, i.e.,  $O^{-1} =O^T$.  Consequently, Lemmas 1 and 2 apply with $T=O$ and $S=O^{-1}_1 $ $O$ both orthogonal.  Equation (\ref{2.4}) holds so $S$ is block diagonal with blocks $S_{ii}$ which are orthogonal and of order $k_i$.

The dimension of the set of orthogonal matrices of order $k_i$ is $k_i (k_i-1)/2$.   Thus the dimension of the set of orthogonal $S$ that commute with $\Lambda$, when the multiplicity of $\Lambda$ is $\left( k_1, k_2, \cdots, k_I\right)$, is 

\noindent{\bf Lemma 5.}
\be
\mbox{real dim} \left\{ S\left| \right. SS^T =I, \:S\Lambda =\Lambda   S , \,\mbox { mult  }\Lambda =\left( k_1, k_2, \cdots, k_I\right)\right\} = \sum\limits^I_{i=1} k_i \left( k_i -1\right) /2.
\label{6.1}
\ee
By using (\ref{6.1}) in the result (\ref{2.2}), we have

\noindent {\bf Theorem  5.}
\be
\mbox{real dim}\left\{ A\left| \right.  A=O\Lambda O^T, \,\mbox{mult }\Lambda =\left(k_1, k_2,\cdots, k_I\right) \right\} =n(n-1)/2 -\sum\limits^I_{i=1} k_i \left( k_i-1\right) /2.
\label{6.2}
\ee
Furthermore, the number of  distinct eigenvalues in $\Lambda$ is $I$.  Therefore by adding $I$ to the dimension in (\ref{6.2}), we get the real dimension of the set of real symmetric  $A$ with multiplicities $k_1, k_2,\cdots, k_I$:

\noindent{\bf Corollary 5.1.}
\be
\mbox{real dim}\left\{ A\left| \right.  A \mbox{ real symm}, \,\mbox{mult}\left( k_1, k_2, \cdots, k_I\right) \right\} =n(n-1) /2 +I -\sum\limits^I_{i=1} k_i \left(k_i-1\right)/2.
\label{6.3}
\ee

Next we subtract the dimension given in (\ref{6.3}) from the dimension of the space of real symmetric matrices of order $n$, which is $n(n+1) /2$, and we state the result as

\noindent {\bf Corollary 5.2.}  In the  space of real symmetric matrices $A$ of order $n$,  the real codimension of the set of matrices having eigenvalue multiplicities $k_1, k_2, \cdots, k_I$ with $k_1 + k_2 +\cdots +k_I=n$, is
\bea
&& \mbox{real codim }\left\{A\left| \right.  A \mbox{ real symm},  \,\mbox{mult}\left( k_1, k_2,\cdots, k_I\right) \right\}\qquad  \qquad  \qquad  \qquad  \qquad  \qquad  \nn\\
&&  \qquad  =
\frac{n(n+1)}{2} -\left[\frac{n(n-1)}{2} \,+\, I -\sum\limits^I_{i=1} \: \frac{k_i (k_i-1)}{2}\right] = \frac{1}{2} \sum\limits^I_{i=1} \left( k_i +2\right) \left( k_i -1\right).
\label{6.4}
\eea
This is the number of conditions which must be  satisfied by the elements of a real symmetric $n$-th order matrix in order for it to have eigenvalues with the specified multiplicities.

From (\ref{6.4}) we can calculate the number of conditions which $A$ must satisfy to have one eigenvalue of multiplicity $k$.  We set $k_1 =k$ and minimize (\ref{6.4}) with respect to the other $k_i$, which requires $k_2 =\cdots=k_I=1$.  Thus we get  

\noindent {\bf Corollary 5.3.}
\be
\mbox{ real codim }\left\{ A\left| \right.  A\mbox{ real symm},  \,\mbox{mult}\left( k, 1,\cdots,1\right)\right\} = \left( k+2\right) \left( k-1\right) /2.
\label{6.5}
\ee
This codimension is independent of $n$.  When $k=2$ it yields the value $2$, obtained by von Neuman and Wigner \cite{vN&W}.

\section{Jordan Canonical Form}
\sce

We shall now use Lemma 2 to determine the dimension of the set of $A$ in $C^{nn}$ which have the Jordan canonical form $J$.  To find the set of $S$ which satisfy (\ref{2.1}), we let $J$ have the block diagonal form
\be
J=\left( J^1, J^2,\cdots,J^K\right), \quad J^i =\lambda_i \, I^i +H^i,
\quad I^i =\left(\begin{array}{ccc}
1& &0\\
&   \ddots& \\
0& &1\end{array}\right), 
\quad H^i \left( \begin{array}{cccc}
\vspace{-10pt} 0 & 1& & 0\\
 \vspace{-10pt}&\ddots&\!\!\!\!\!\ddots&\\
 \vspace{-10pt} &&\ddots&1\\
 0&&&0 \end{array}\right).
\label{7.1}
\ee
Here $J^i$, $I^i$ and $H^i$ are square matrices  of order $k_i$.  We write $S= \left(S^{ij}\right)$ where the block $S^{ij}$ is $k_i$ by $k_j$.  Then  $\left( SJ\right)^{ij} =S^{ij} J^j$ and $\left( J S\right)^{ij} =J^i S^{ij}$, so (\ref{2.1})  becomes
\be
S^{ij} J^j =J^i S^{ij}.
\label{7.2}
\ee

Upon using the definition (\ref{7.1})  of $J^i$ in (\ref{7.2}) we get
\be
S^{ij} \lambda_j +S^{ij} H^j =\lambda_i S^{ij} + H^i S^{ij}.
\label{7.3}
\ee
Next we use the definition of $H^i$ in (\ref{7.3}), and we write the $st$ element of $S^{ij}$ as $S^{ij}_{st}$.  In this way we obtain
\bea
S^{ij}_{st} \left( \lambda_j -\lambda_i\right) &=&-\sum\limits_q S^{ij}_{sq} H^j_{qt} + \sum\limits_r H^i_{sr} S^{ij}_{rt} =-\sum\limits_qS^{ij}_{sq} \,\delta_{q, t-1} + \sum\limits_r \delta_{s, r-1} S^{ij}_{rt}\nn\\
&=& -S^{ij}_{s, t-1} + S^{ij}_{s+1, t}.
\label{7.4}
\eea

For $\lambda_i\not= \lambda_j$, (\ref{7.4}) determines $S^{ij}_{st}$ in terms of  $S^{ij}_{s, t-1}$ and $S^{ij}_{s+1, t}$.  Repeated use of (\ref{7.4}) leads either  to $t-1 =0$ or to $s+1 =k_i+1$.  But $S^{ij}_{so}=0$ and $S^{ij}_{k_i+1, t} =0$ so
\be
S^{ij}_{st} =0 \qquad \mbox{ for } \lambda_i \not= \lambda_j.
\label{7.5}
\ee
If $\lambda_i =\lambda_j$ then (\ref{7.4}) yields $\displaystyle{S^{ij}_{s, t-1} =S^{ij}_{s+1, t}}$, so $S^{ij}$ is a Toeplitz matrix, i.e., it is constant on lines parallel to the main diagonal.  The entries are zero where $t<s+\left[k_j -k_i\right]_+$:
\be
S^{ij}_{st} =0 \mbox{ if } t<s+\left[ k_j -k_i\right]_+ \mbox{ when }\lambda^i =\lambda^j.
\label{7.6}
\ee
Thus for $\lambda_i =\lambda_j$, $S^{ij}$ has the first form in (\ref{7.7}) for $k_i\leq k_j$ and the second form in (\ref{7.7}) for $k_i \geq k_j$:

\begin{figure}[h]
\unitlength=1in
\begin{center}
\psfragscanon
\includegraphics[width=4in]{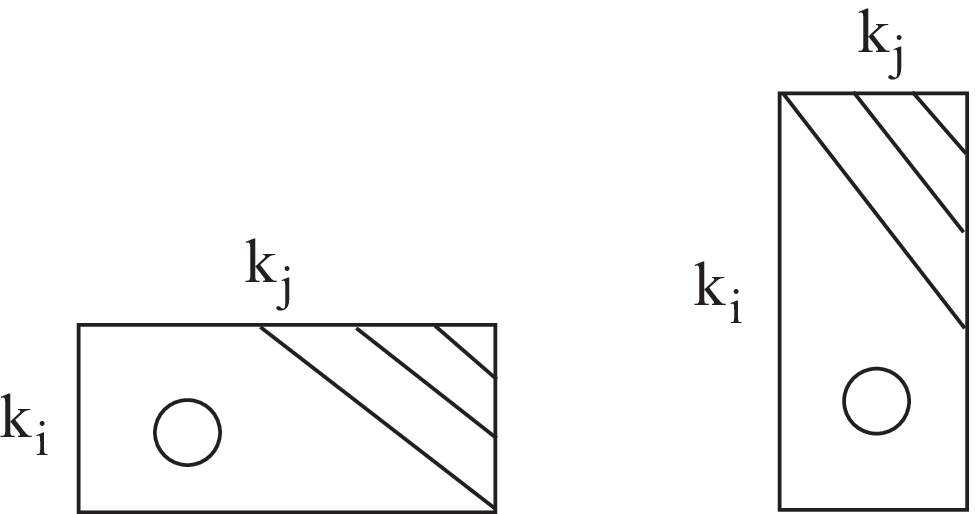}
\end{center}
\end{figure}

\vskip-1.5in 
\be
\label{7.7}
\ee
\vskip1in 
\noindent This form of the $S^{ij}$ was found by Gantmacher  \cite[p.\ 221]{Gantmacher}.

The complex dimension of the set of $S^{ij}$ satisfying (\ref{7.6}) is 
\be
\mbox{ complex dim } \left\{ S^{ij} \mbox{ satisfying  (\ref{7.6})}\right\} =\left\{ \begin{array}{lcl}
k_j -\left[ k_j -k_i\right]_+ &=& k_i \mbox{ if } k_i \leq k_j\\
k_j -0 &=&k_j \mbox{ if } k_i > k_j\end{array} \right\}
=\min \left( k_i , k_j\right).
\label{7.8}
\ee
Thus the dimension of the set of all $S$ satisfying (\ref{2.1}) is
\be
\mbox{complex dim } \left\{ S\left| \right.  SJ=JS\right\} =\sum\limits_{i,j} \mbox{ dim } \left\{ S^{ij}\right\} =\mathop{\sum\limits_{i,j}}\limits_{\lambda_i \not= \lambda_j} \; 0 + \mathop{\sum\limits_{i,j}}\limits_{\lambda_i = \lambda_j} \min \left( k_i , k_j\right).
\label{7.9}
\ee

We can write the sum in (\ref{7.9}) more explicitly  by specifying that there are $p$ distinct eigenvalues $\lambda_a$, and that $\lambda_a$ occurs in $N_a$ blocks.  We write their orders $k_{a1} \geq k_{a2} \geq\cdots\geq k_{aN_{a}}$, 
and set $k_{aj}=0$ for $j>N_a$.  Then $k_{a1} +\cdots+k_{aN_a}=n_a $, where $n_a$ is the  multiplicity of $\lambda_a$.  Now
 we set $\lambda_1=\lam_j =\lam_a$ in (\ref{7.9}), replace $k_i$ and $k_j$ by $k_{ai}$ and $k_{aj}$, and sum over $a$.  Thus we get
\bea
\mbox{complex dim } \left\{ S\left|\right. SJ = JS\right\}^i   
&=& \mathop{\sum\limits_{i,j}}\limits_{\lam_i = \lam_j}  \min \left( k_i, k_j\right)\nn\\
&=&
\sum\limits^p_{a=1}\sum\limits^{Na}_{i=1}\sum\limits^{Na}_{j=1} \min \left( k_{ai}, k_{aj}\right) = \sum\limits^p_{a=1}\sum\limits^{Na}_{j=1} \left( 2j -1\right) k_{aj}\nn\\
&=& \sum\limits^{N^\ast}_{j=1} \left( 2j-1\right) \sum\limits^p_{a=1} k_{aj}.
\label{7.10}
\eea
Here $N^\ast =\max\limits_a \, N_a$ is the maximum number of blocks with a given eigenvalue.  The final interchange of the order of summation is valid because $k_{aj} =0$ for $j>N_a$.

The sum of the $k_{aj}$ in (\ref{7.10}) is just the degree $m_j$ of the $j$-th invariant polynomial of $J$, and therefore also of any $A$ which is similar to $J$. Therefore (\ref{7.10}) can be rewritten in the following form (Gantmacher \cite[p.\ 222, Theorem 2]{Gantmacher}):
\be
\mbox{complex dim } \left\{ S\left| \right.  SJ=JS\right\} =\sum\limits^{N^\ast}_{j=1} \left( 2j-1\right) m_j.
\label{7.11}
\ee

Upon using (\ref{7.11}) in (\ref{2.2}), and recalling that $\mbox{complex dim}\left\{ T\right\}=n^2$, we have 

\noindent{\bf Theorem 6.}  The complex dimension of the  set of $A\in C^{nn}$ with the Jordan form  $J$ is 
\be
\mbox{ complex dim } \left\{ A\left| \right. A=TJT^{-1}\right\} =n^2
- \sum\limits^{N^\ast}_{j=1} \left( 2j -1\right) \sum\limits^p_{a=1} k_{aj} = n^2 -\sum\limits^{N^\ast}_{j=1} \left( 2j-1\right) m_j.
\label{7.12}
\ee
The complex dimension of the set of $p$ distinct eigenvalues is $p$.  Upon adding $p$ to (\ref{7.12}) we get the complex dimension of the set of $A$   with the structure of $J$, but with any $p$ distinct eigenvalues.  We can state this result as follows:

\noindent{\bf Theorem 7.}  Let $A\in C^{nn}$ have $p$ distinct eigenvalues $\lambda_a$ with multiplicities $n_a$, $a=1,\cdots,p$ such that  $n_1+n_2+\cdots+n_p=n$.  Suppose that $\lambda_a$ occurs in $N_a$ blocks of the Jordan canonical form of $A$, with the blocks of orders $k_{a1} \geq k_{a2}\cdots \geq k_{aN_a}$, with $k_{a1} + k_{a2}+\cdots+k_{aN_a}=n_a$.  The complex dimension of the set of $A$ with these properties, and any values of the $\lambda_a$, is
\bea
\mbox{ complex dim } \left\{ A\left| \right. n_a, N_a, k_{aj}\right\}
&=& n^2 -\sum\limits^{N^\ast}_{j=1} \left( 2j -1\right) \sum\limits^p_{a=1} \; k_{aj} +p\nn\\
&=& n^2 -\sum\limits^{N^\ast}_{j=1} \left( 2j-1\right) m_j +p.
\label{7.13}
\eea

When all the eigenvalues of $A$ are simple, then $p=n$, $n_a =N_a =N^\ast =k_{aj} =1$ and (\ref{7.13}) yields $\mbox{ complex dim }\left\{ A \right\} =n^2$.  When $A$ has one $n$-fold eigenvalue and just one block in its Jordan form, then $p=1$, $n_1=n$, $N^\ast = N_1 =1$, $k_{11} =n$ and (\ref{7.13}) yields
\be
\mbox{ complex dim } \left\{ A \left| \right. p=1, \: n_1=n,\:  N_1=1, \: k_{11} =n \right\} =n^2-n+1.
\label{7.14}
\ee
When $A$ has one $n$-fold eigenvalue and $n$ blocks in its Jordan form, then $p=1$, $n_1=n$, $ N^\ast = N_1 =n$, $k_{11}=1$ and (\ref{7.13}) yields 
\be
\mbox{ complex dim } \left\{ A \left| \right. p=1, \:n_1=n,\: N_1=n, \: k_{11} =1 \right\} =1.
\label{7.15}
\ee

\section{Singular value decomposition}
\sce
Every $A\in R^{nm}$, the space of real $n$ by $m$ matrices, has a singular value decomposition
\be
A=U\Sigma V^T.
\label{8.1}
\ee
Here $U\in R^{nn}$ and $V\in R^{mm}$ are orthogonal.  $\Sigma\in R^{nm}$ is diagonal with $r\leq \min (n, m)$ nonzero elements, on the main diagonal, starting at the upper left corner.  They are all positive, and are called singular values of $A$.  We seek the real dimension of the set of $A$ with $J\leq r$ distinct  singular values $\sigma_1> \sigma_2>\cdots>\sigma_J >0$, with respective multiplicities $k_1, k_2,\cdots,k_J$.  Then $r=k_1+k_2+\cdots+k_J$.

Suppose that in addition to (\ref{8.1}), $A$ can be written as $A= U_1 \Sigma V^T_1$ where $U_1$ and $V_1$ are also orthogonal.  Then we    equate these two expressions for $A$ to get $U\Sigma V^T=U_1 \Sigma V^T_1$.  From this we conclude that $U^T_1 U\Sigma =\Sigma V^T_1 V$, which we rewrite as
\be
Q\Sigma =\Sigma P.
\label{8.2}
\ee
Here $Q=U^T_1U\in R^{nn}$ and $P=V^T_1 V\in R^{mm}$ are both orthogonal, being   products  of orthogonal factors.  Then $U_1 =UQ^T$ and $V_1 =VP^T$ give the same $A$ as do $U$ and $V$ when used in (\ref{8.1}).

From (\ref{8.1}) and this conclusion we obtain the following lemma:

\noindent {\bf Lemma 5.}  The real dimension of the set of $A\in R^{nm}$ with the singular value matrix $\Sigma$ is
\bea
&&\mbox{\hskip-0.5in }\mbox{ dim } \left\{ A\left| \right.  A\in R^{nm}, \quad A=U\Sigma V^T\right\}\nn\\
&&  = \dim \left\{ U\left| \right. U\in R^{nn},\quad UU^T=I\right\} + \mbox{ dim }\left\{ V\left| \right. V\in R^{mm},\quad VV^T=I\right\}\nn\\
&&- \mbox{ dim } \left\{ Q, P\left| \right. Q\in R^{nn},\quad QQ^T =I, \quad P\in R^{mm}, \quad PP^T=I, \quad Q\Sigma =\Sigma P\right\}.
\label{8.3}
\eea
The real dimensions of the sets of orthogonal $U$ and $V$ are $n(n-1)/2$ and $m(m-1)/2$ respectively.  Therefore to use the lemma we shall determine the set of pairs $Q,P$ and its dimension.

We begin by writing the elements $\Sigma_{ij}$ of $\Sigma$ in the form
\begin{displaymath}
\Sigma_{ii} =\sigma^\pr_i, \quad i=1,\cdots, r; \quad \Sigma_{ij} =0 \: \mbox{ if } i\not= j \mbox{ or if }i=j>r.
\end{displaymath}
Here the first  $k_1$ of the $\sigma^\pr_i$ are equal to $\sigma_1$, the next $k_2$ of the $\sigma^\pr_i$ equal  $\sigma_2$, etc.   Then we find that 
\bea
\left( Q\Sigma\right)_{ik} &=& \left\{ \begin{array}{ll} Q_{ik} \sigma^\pr_k, & k\leq r\\
0,& k>r\end{array}\right.
\label{8.4}\\
\left( \Sigma P\right)_{ik} &=& \left\{ \begin{array}{ll} \sigma^\pr_i P_{ik},  & i\leq r\\
0,&i>r.
\end{array}\right.
\label{8.5}
\eea
Now we use (\ref{8.4}) and (\ref{8.5}) in (\ref{8.2}) to find
\begin{subequations}
\bea
Q_{ik} \sigma^\pr_k = \sigma^\pr_i P_{ik}, && i\leq r, \quad k\leq r\label{8.6a}\\
Q_{ik} \sigma^\pr_k =0,& &  i>r,\quad  k\leq r\label{8.6b}\\
0= \sigma^\pr_i P_{ik},  & &  i\leq r, \quad k>r.
\label{8.6c}
\eea
\end{subequations}
Since the $\sigma^\pr_i\not= 0$, we obtain from (8.6)
\begin{subequations}
\bea
Q_{ik} =\frac{\sigma^\pr_i}{\sigma^\pr_k}\, P_{ik},  & & i,k\leq r\label{8.7a}\\
Q_{ik} =0, && i>r, \quad k\leq r\label{8.7b}\\
P_{ik} =0, && i\leq r,\quad  k>r.
\label{8.7c}
\eea
\end{subequations}

From (8.7) and the orthogonality of $Q$ and $P$, one can prove

\noindent{\bf Theorem 8.}  The orthogonal matrices $Q$ and $P$ satisfying (\ref{8.2}) are  block diagonal, each with $J+1$   square blocks $Q_1,\cdots,Q_{J+1}$ and $P_{1},\cdots,P_{J+1}$.  Each block is an orthogonal matrix.  For $j=1, \cdots, J$ the blocks $Q_j$ and $P_j$ are of  order $k_j$ and they are equal:  $Q_j=P_j$, $j=1,\cdots,J$.  $Q_{J+1}$ is of order $n-r$ and $P_{J+1}$ is of order $m-r$.

This theorem is given by Horn and Johnson \cite[p.\ 147, (3.1.1$^\pr$)]{horn&johnson}.

From Theorem 8 we obtain

\noindent {\bf Corollary 8.1.}
\bea
&&\mbox{ dim } \left\{ Q,P \left| \right.  Q\in R^{nn},\quad QQ^T=I, \quad P\in R^{mm}, \quad P P^T =I, \quad Q\Sigma =\Sigma P\right\}\nn\\
&&\qquad  =\sum\limits^J_{j=1} \frac{k_j (k_j-1)}{2} +\frac{(n-r) (n-r-1)}{2} + \frac{(m-r)(m-r-1)}{2}.
\label{8.8}
\eea

We now use (\ref{8.8}) in (\ref{8.3}), and simplify the result to get

\noindent {\bf Theorem 9.}  The real dimension of the set of $A\in R^{nm}$ with the singular values $\sigma_1,\cdots,\sigma_J$ having multiplicities $k_1,\cdots,k_J$ is
\bea
&&\!\!\!\!\!\!\!\!\!\! \mbox{ dim } \left\{ A\left| \right. A\in R^{nm}, \quad A =U\Sigma V^T\right\} \nn\\
&& =\frac{1}{2} n (n-1) +\frac{1}{2}m (m-1) -\frac{1}{2} (n-r) (n-r-1) -\frac{1}{2}(m-r)(m-r-1) -\frac{1}{2} \sum\limits^J_{j=1} \; k_j (k_j-1)\nn\\
&&= (n+m-r) r-r-\frac{1}{2} \sum\limits^J_{j=1} k_j (k_j-1).
\label{8.9}
\eea
By adding to (\ref{8.9}) the dimension $J$ of the set of singular values, we get the dimension of the set of $A$ having singular values with the specified multiplicities $k_1,\cdots,k_J$:

\noindent{\bf Corollary 9.1.}
\be
\mbox{ dim } \left\{ A\left| \right.  A\in R^{nm} , \,\mbox{mult}\left( k_1,\cdots,k_J\right)\right\} =(n+m-r) r-r -\frac{1}{2} \sum\limits^J_{j=1} k_j (k_j-1) +J.
\label{8.10}
\ee
The codimension in $R^{nm}$ is just $nm$ minus the right side of (\ref{8.10}).

As an example, if all the $\sigma_j$ are simple then $J=r$ and each $k_j=1$, so (\ref{8.10}) yields 
\be
\mbox{ dim } \left\{ A\left| \right. A\in R^{nm},  \,\mbox{mult}\left( 1,\cdots,1\right)\right\} =(n+m-r)r.
\label{8.11}
\ee
This is just the dimension of the set of $A\in R^{nm}$ of rank $r$.  The codimension of the set of $A$ in Corollary   9.1 with respect to this set of matrices is just the right side of (\ref{8.11}) minus the right side of (\ref{8.10}):
\bea
&&\mbox{ codim (with respect to rank $r$ matrices)} \; \left\{ A\left| \right.  A\in R^{nm},  \,\mbox{mult}\left( k_1,\cdots,k_J\right)\right\}\nn\\
&&\qquad  \qquad  =\frac{1}{2} \sum\limits^J_{j=1} k_j (k_j-1) + r-J.
\label{8.12}
\eea

\newpage

\begin{center}
\begin{tabular}{rlcc}
\multicolumn{4}{c} {\bf Table 1}\\
&\raisebox{-1.5ex}[0pt]{Set of} &\raisebox{-1.5ex}[0pt]{Complex }&\raisebox{-1.5ex}[0pt]{ Real}\\
&\raisebox{1.5ex}[0pt]{Matrices in $C^{nn}$ }&   \raisebox{1.5ex}[0pt]{dimension} &\raisebox{1.5ex}[0pt]{ dimension}\\
\hline
1. & All & $n^2$ & $2n^2$\\
2. & Invertible & $n^2$ & $2n^2$\\
3. & Singular & $n^2-1$ & $2(n^2-1)$\\
4. & Diagonalizable & $n^2$ & $2n^2$\\
5. & Normal & $n(n+1)/2$ & $n(n+1)$\\
6. & Hermitian & --- & $n^2$\\
7. & Unitary & --- & $n^2$ \\
8. & Symmetric & $n(n+1)/2$ & $n(n+1)$\\
9. & Real Symmetric & --- & $n(n+1)/2$\\
10. & Antisymmetric & $n(n-1)/2$ & $n(n-1)$\\
11. & Real Antisymmetric & --- & $n(n-1)/2$\\
12. & Orthogonal & --- & $n(n-1)/2$\\
13.& Matrices in $C^{nn}$   &&\\
&\raisebox{1.5ex}[0pt]{rank $r\leq \min (n, m)$} & \raisebox{1.5ex}[0pt]{$(m+n-r)r$ }& \raisebox{1.5ex}[0pt]{$2(m+n-r)r$}
\end{tabular}
\end{center}

\vskip.5in 

\begin{center}
\begin{tabular}{rlcc}
\multicolumn{4}{c}{\bf Table 2}\\
\multicolumn{2}{c}{\raisebox{-1.5ex}[0pt]{\hskip.25in Set of matrices in $C^{nn}$ }}&\raisebox{-1.5ex}[0pt]{Complex }& \raisebox{-1.5ex}[0pt]{Real} \\
\multicolumn{2}{c}{\raisebox{1.5ex}[0pt]{mult $k_1,\cdots,k_I$ }}&\raisebox{1.5ex}[0pt]{dimension} & \raisebox{1.5ex}[0pt]{dimension}\\
\hline
1. & Diagonalizable & $n^2- \sum\limits^I_{i=1} (k^2_i -1)$ & $ 2\left[ n^2 -\sum\limits^I_{i=1} (k^2_i -1)\right]$\\
2. & Normal & --- & $n^2 +I -\sum\limits^I_{i=1} \left( k^2_i -1\right)$\\
3. & Hermitian & --- & $n^2 -\sum\limits^I_{i=1} \left( k^2_i -1\right)$\\
4. & Unitary & --- & $n^2 - \sum\limits^I_{i=1} \left( k^2_i -1\right)$\\
5. &     Real Symmetric & --- & $\frac{n(n-1)}{2} +I-\frac{1}{2} \sum\limits^I_{i=1} k_i 
\left( k_i-1\right)$\\
6. & Matrices in $R^{nm}$ &--- & $(n+m -r) r-r+I -\frac{1}{2} \sum\limits^I_{i=1} k_i 
 \left( k_i -1\right)$\\
& \raisebox{1.5ex}[0pt]{mult $k_1,\cdots,k_I$}&&\\
&\raisebox{1.5ex}[0pt]{ $\sum\limits^I_{i=1}\, k_i =r$}&&\\
7. & with normal form $J$ & ---& $n^2 -\sum\limits^{N^\ast}_{j=1} \left( 2j-1\right) \sum\limits^p_{a=1} k_{aj} + p$\\
&\raisebox{1.5ex}[0pt]{(Theorem 7)} &&\\
8. & $A\in R^{nm}$ \, with singular value &---& $(n+m-r) r-r -\frac{1}{2} \sum\limits^J_{j=1} \, k_j \left( k_j -1\right) +J$\\
&\raisebox{1.5ex}[0pt]{multiplicities $k_1,\cdots,k_J$} &&\\
&\raisebox{1.5ex}[0pt]{(Corollary 9.1)}&&
\end{tabular}
\end{center}

\end{document}